\documentclass[11pt]{article}
\usepackage{epic,latexsym,amssymb}
\usepackage{amsfonts}
\usepackage{amscd}
\usepackage{amsmath}
\usepackage{graphicx}
\usepackage{color}
\usepackage{caption}
\usepackage{tikz}

\usepackage{float}

\newtheorem{thm}{Theorem}

\newtheorem{ob}[thm]{Observation}

\newtheorem{lem}[thm]{Lemma}

\newcommand{\diam}{{\rm diam}}

\newcommand{\cP}{\mathcal{P}}

\newcommand{\barS}{\overline{S}}

\newcommand{\proof}{\noindent\textbf{Proof. }}

\newcommand{\qed}{$\Box$}

\newcommand{\1}{\vspace{0.1cm}}
\newcommand{\2}{\vspace{0.2cm}}

\newcommand{\QEDmark}{\mbox{\textsc{qed}}}
\newcommand{\proofStarter}[1]{\textsc{#1} }

\def\vertex(#1){\put(#1){\circle*{2}}}
\def\vertexo(#1){\put(#1){\circle{2}}}
\def\vert(#1){\put(#1){\circle*{1.5}}}
\def\verto(#1){\put(#1){\circle{1.5}}}
\def\lab(#1)#2{\put(#1){\makebox(0,0)[c]{#2}}}
\definecolor{DarkGreen}{rgb}{0.2, 0.6, 0.3}

\definecolor{electricindigo}{rgb}{0.44, 0.0, 1.0}

\newcommand{\gtd}{\gamma_2}

\let\oldenumerate\enumerate
\renewcommand{\enumerate}{
  \oldenumerate
  \setlength{\itemsep}{0.5pt}
  \setlength{\parskip}{0pt}
  \setlength{\parsep}{0pt}
}

\begin{document}

\title{On the Total Forcing Number of a Graph}
\author{$^{1,2}$Randy Davila and $^{1}$Michael A. Henning\\
\\
$^1$Department of Pure and Applied Mathematics\\
University of Johannesburg \\
Auckland Park 2006, South Africa \\
\\
$^2$Department of Mathematics \\
Texas State University \\
San Marcos, TX 78666, USA \\
\small {\tt Email: rrd32@txstate.edu}
}

\date{}
\maketitle

\begin{abstract}
A total forcing set in a graph $G$ is a forcing set (zero forcing set) in $G$ which induces a subgraph without isolated vertices. Total forcing sets were introduced and first studied by Davila~\cite{Davila}. The total forcing number of $G$, denoted $F_t(G)$ is the minimum cardinality of a total forcing set in $G$. We study basic properties of $F_t(G)$, relate $F_t(G)$ to various domination parameters, and establish $NP$-completeness of the associated decision problem for $F_t(G)$. Our main contribution is to prove that if $G$ is a connected graph of order~$n \ge 3$ with maximum degree~$\Delta$, then $F_t(G) \le ( \frac{\Delta}{\Delta +1} ) n$, with equality if and only if $G$ is a complete graph $K_{\Delta + 1}$, or a star $K_{1,\Delta}$.
\end{abstract}


\section{Introduction}

In recent years, dynamic colorings of the vertices in a graph has gained much attention. Indeed, forcing sets \cite{AIM-Workshop, zf_np, Davila Kenter, Genter1, Genter2, Edholm, Hogben16, zf_np2}, $k$-forcing sets \cite{k-Forcing, Dynamic Forcing}, connected forcing sets \cite{CF Complexity, Brimkov Davila, DaHeMaPe16}, and power dominating sets \cite{powerdom3, powerdom2}, have seen a wide verity of application and interesting relationships to other well studied graph properties. These aforementioned sets all share the common property that they may be defined as graph colorings that change during discrete time intervals. Of these dynamic colorings, we highlight that the most prominent is that of \emph{forcing} (\emph{zero forcing}), and the associated graph invariant known as the \emph{forcing number} (\emph{zero forcing number}). This paper continues the study of forcing in graphs by way of restricting the structure of forcing sets as an induced subgraph.

Let $G = (V,E)$ be a graph with vertex set $V = V(G)$ and edge set $E = E(G)$. The \emph{forcing process} is defined as follows: Let $S\subseteq V$ be a set of initially ``colored" vertices, all other vertices are said to be ``non-colored". A vertex contained in $S$ is said to be $S$-colored, while a vertex not in $S$ is said to be $S$-uncolored. At each time step, if a colored vertex has exactly one non-colored neighbor, then this colored vertex \emph{forces} its non-colored neighbor to become colored. If $v$ is such a colored vertex, we say that $v$ is a \emph{forcing vertex}. We say that $S$ is a \emph{forcing set}, if by iteratively applying the forcing process, all of $V$ becomes colored. We call such a set $S$ an $S$-forcing set. In addition, if $S$ is an $S$-forcing set in G and $v$ is a $S$-colored vertex that forces a new vertex to be colored, then $v$ is an \emph{$S$-forcing vertex}. The cardinality of a minimum forcing set in $G$ is the \emph{forcing number} of $G$, denoted $F(G)$. If $S$ is a forcing set which  also induces a connected subgraph, then $S$ is a \emph{connected forcing set}. The cardinality of a minimum connected forcing set in $G$ is the \emph{connected forcing number} of $G$, denoted $F_c(G)$.

For graphs in general, it is known that computation of both $F(G)$ and $F_c(G)$ lie in the class of $NP$-hard decision problems, see \cite{zf_np, zf_np2} and \cite{CF Complexity}, respectively. Moreover, $F(G)$ and $F_c(G)$ have been related to many well studied graph properties such as minimum rank, independence, and domination, see for example \cite{AIM-Workshop, k-Forcing, DaHeMaPe16}. For more on forcing and connected forcing, we refer the reader to \cite{k-Forcing, CF Complexity, Brimkov Davila, Dynamic Forcing, Davila, Davila Henning, DaHeMaPe16, Davila Kenter, Genter1, Genter2}.

In this paper we study a variant of forcing. Namely, if $S\subseteq V$ is a forcing set of $G$ that induces a subgraph without isolated vertices, then $S$ is a \emph{total forcing set}, abbreviated TF-\emph{set}, of $G$. The concept of a total forcing set was first introduced and studied by the Davila in~\cite{Davila}. The minimum cardinality of a TF-set in $G$ is the \emph{total forcing number} of $G$, denoted $F_t(G)$. Minimum cardinality TF-sets in $G$ are called $F_t(G)$-sets.

We proceed as follows. In the next section concepts used throughout the paper are introduced and known facts and results needed are recalled. In Section~\ref{S:properties} we provide fundamental properties of total forcing sets in graphs. Relationships between total forcing and various domination parameters is explored in Section~\ref{S:TotalFvsDom}. In Section~\ref{S:vertexRem}, we study the effect on the total forcing number when the vertex whose removal creates no isolates is deleted from a graph. The computational complexity of the total forcing number is discussed in Section~\ref{S:NP} where it is shown that the decision problem associated with total forcing is also $NP$-complete. We close in Section~\ref{S:bds} with an upper bound on the total forcing number of a graph with minimum degree at least two in terms of the order and maximum degree of the graph.

\section{Definition and Known Results}

For notation and graph terminology, we will typically follow \cite{MHAYbookTD}. Throughout this paper, all graphs will be considered undirected, simple and finite. Specifically, let $G$ be a graph with vertex set $V(G)$ and edge set $E(G)$ of order $n = |V(G)|$ and size $m = |E(G)|$.
Two vertices $v$ and $w$ are \emph{neighbors} in $G$ if they are adjacent; that is, if $vw \in E(G)$. The \emph{open neighborhood} of a vertex $v$ in $G$ is the set of neighbors of $v$, denoted $N_G(v)$, whereas the \emph{closed neighborhood} of $v$ is $N_G[v] = N_G(v) \cup \{v\}$. The \emph{open neighborhood} of a set $S\subseteq V(G)$ is the set of all neighbors of vertices in $S$, denoted $N_G(S)$, whereas the \emph{closed neighborhood} of $S$ is $N_G[S] = N_G(S) \cup S$. The \emph{degree} of a vertex $v$ in $G$, is denoted $d_G(v) = |N_G(v)|$. The minimum and maximum degree of $G$ are denoted by $\delta(G)$ and $\Delta(G)$, respectively. For a subset $S \subseteq V(G)$, the \emph{degree of $v$ in $S$}, denoted $d_S(v)$, is the number of vertices in $S$ adjacent to~$v$; that is, $d_S(v) = |N(v) \cap S|$. In particular, $d_G(v) = d_{V(G)}(v)$.
If the graph $G$ is clear from the context, we simply write $V$, $E$, $n$, $m$, $d(v)$, $N(v)$, $N(S)$, $\delta$ and $\Delta$ rather than $V(G)$, $E(G)$, $n(G)$, $m(G)$, $d_G(v)$, $N_G(v)$, $N_G(S)$, $\delta(G)$ and $\Delta(G)$, respectively.

The distance between two vertices $v,w\in V$ is the length of a shortest $(v,w)$-path in $G$, and is denoted by $d_G(v,w)$. If no $(v,w)$-path exists in $G$, then we define $d_G(v,w) =\infty$. The maximum distance among all pairs of vertices of $G$ is the \emph{diameter} of $G$, denoted by $\diam(G)$. The length of a shortest cycle in a graph $G$ (containing a cycle) is the \emph{girth} of $G$, denoted by $g = g(G)$. For a set of vertices $S \subseteq V$, the subgraph induced by $S$ is denoted by $G[S]$. The subgraph obtained from $G$ by deleting all vertices in $S$ and all edges incident with vertices in $S$ is denoted by $G - S$. If $S = \{v\}$, we simply write $G - v$ rather than $G - S$. We will denote the path, cycle, and complete graph on $n$ vertices by $P_n$, $C_n$, and $K_n$, respectively. A \emph{leaf} of $G$ is a vertex of degree~$1$ in $G$.

A \emph{packing} in a graph is a set of vertices that are pairwise at distance at least~$3$ apart; that is, if $P$ is a packing in a graph $G$, and $u$ and $v$ are distinct vertices of $P$, then $d_G(u,v) \ge 3$. We note that if $P$ is a packing, then the closed neighborhoods, $N_G[v]$, of the vertices $v$ in $P$ are pairwise vertex disjoint. A \emph{perfect packing} (also called a \emph{perfect dominating set} in the literature) is a packing that dominates the graph; that is, if $P$ is a perfect packing, then the closed neighborhoods, $N_G[v]$, of the vertices $v$ in $P$ partition $V(G)$.

We use the standard notation $[k] = \{1,\ldots,k\}$.

\medskip
\noindent\textbf{Domination in Graphs.}  A set of vertices $S\subseteq V$ is a \emph{dominating set}, if every vertex not in $S$ has a neighbor in $S$. The minimum cardinality of a dominating set in $G$ is the \emph{domination number} of $G$, denoted by $\gamma(G)$. If $S\subseteq V$ has the property that every vertex in $G$ has a neighbor in $S$, then $S$ is a \emph{total dominating set}. The minimum cardinality of a total dominating set in $G$ is the \emph{total domination number} of $G$, denoted by $\gamma_t(G)$. If $S\subseteq V$ is a dominating set with the additional property that $S$ induces a connected subgraph, then $S$ is a \emph{connected dominating set}. The minimum cardinality of a connected dominating set in $G$ is the \emph{connected domination number} of $G$, denoted $\gamma_c(G)$. A set $S$ of vertices in a graph $G$ is a (distance) $2$-\emph{dominating set} if every vertex not in $S$ is within distance~$2$ from some vertex of $S$. The (distance) $2$-\emph{domination number}, written $\gtd(G)$, is the minimum cardinality of a $2$-dominating set in $G$.

Domination and its variants are heavily studied in graph theory and we refer the reader to the monographs \cite{DomBook1, DomBook2, MHAYbookTD} which detail and survey many results on the topic. A survey on distance domination in graphs can be found in~\cite{He98}.

Serving as a dynamic approach to domination, the power domination process is defined as follows: For a given set of vertices $S\subseteq V$, the sets $(\cP_G^i(S))_{i\ge 0}$ of vertices \emph{monitored} by $S$ at the $i$-th step are defined recursively by, \\[-22pt]
\begin{enumerate}
\item $\cP_G^0(S) = N_G[S]$, and
\item $\cP_G^{i+1}(S) = \cup\{N_G[v] \colon v\in \cP_G^i(S) \colon |N_G[v] \setminus \cP_G^i(S)| \le 1\}.$
\end{enumerate}

If $\cP_G^{i_0}(S) = \cP_G^{i_0+1}(S)$, for some $i_0$, then $\cP_G^j(S) = \cP_G^{i_0}(S)$, for all $j\ge i_0$. We define $\cP_G^{\infty}(S) = \cP_G^{i_0}(S)$. If $\cP_G^{\infty}(S) = V$, then $S$ is a \emph{power dominating set} of $G$. The minimum cardinality of a power dominating set in $G$ is the \emph{power domination number} of $G$, denoted by $\gamma_P(G)$.

Since every dominating set is also a power dominating set, and since every total dominating set is a dominating set, and since every connected dominating set is a total dominating set (provided $\gamma_c(G)\ge 2$), we make note of the following chain of inequalities.

\begin{ob}\label{domOb}
If $G$ is a connected graph of order $n \ge 3$ and $\gamma_c(G)\ge 2$, then
\begin{equation*}
\gamma_P(G) \le \gamma(G) \le \gamma_t(G) \le \gamma_c(G).
\end{equation*}
\end{ob}

\medskip
\noindent\textbf{Known Results on Forcing.} Recalling that computation of forcing for a general graph is $NP$-hard \cite{CF Complexity, zf_np, zf_np2}, we remark that finding computationally efficient lower and upper bounds for $F(G)$ and $F_c(G)$ has been of particular interest, see for example \cite{k-Forcing, CF Complexity, Brimkov Davila, DaHeMaPe16}. In relation to this paper, we highlight that when considering efficient bounds for $F(G)$ and $F_c(G)$, both the minimum degree and maximum degree play an important role.

In the introductory paper on forcing in graphs, which first appeared due to a workshop on the minimum rank of a graph (AIM-Group) \cite{AIM-Workshop}, it was shown that the minimum degree bounds $F(G)$ from below. That is, for any graph $G$ with minimum degree $\delta$, $F(G)\ge \delta$. This minimum degree lower bound has recently been significantly improved when the graph in question has restrictions on its girth and minimum degree. In particular, if $G$ is a graph has minimum degree $\delta \ge 2$, Genter, Penso, Rautenbach, and Souzab \cite{Genter1} and Gentner and Rautenbach \cite{Genter2} proved $F(G) \ge \delta + (\delta -2)(g - 3)$, whenever $g\le 6$, and also whenever $g\le 10$ by Davila and Henning \cite{Davila Henning}. Using the techniques presented in \cite{Davila Henning}, Davila, Malinowski, and Stephen \cite{girthTheorem} recently resolved this inequality in the affirmative for graphs with arbitrary girth.

As shown in \cite{k-Forcing}, Amos, Caro, Davila, and Pepper, provided the first known upper bound on $F(G)$ in terms of maximum degree and order. In particular, if $G$ is an isolate-free graph of order $n$ and maximum degree $\Delta$, they proved $F(G) \le ( \frac{\Delta}{\Delta+1} ) n$, and $F(G) \le \frac{(\Delta -2)n+2}{\Delta -1}$, whenever the added restrictions that $G$ is connected and $\Delta \ge 2$ are imposed. Moreover, they also showed that the forcing number of a graph is related to the connected domination number of $G$ by way of the inequality $F(G) \le n - \gamma_c(G)$. Improving on these aforementioned maximum degree and order upper bounds, Caro and Pepper \cite{Dynamic Forcing} gave a greedy algorithm which resulted in $F(G) \le \frac{(\Delta -2)n -(\Delta -\delta)+2}{\Delta -1}$, whenever $G$ is connected and $\Delta\ge 2$. We remark that an analogous upper bound for $F_c(G)$ does not yet exist. Indeed, Davila, Henning, Magnant, and Pepper \cite{DaHeMaPe16} asked whether or not there exists a function of order and maximum degree which bounds $F_c(G)$ from above.


The main aim of this paper is to introduce and study the total forcing number of a graph. In particular, we establish fundamental properties of $F_t(G)$, establish relationships between $F_t(G)$ and domination, prove $NP$-completeness of the total forcing decision problem, and provide a greedy algorithm for total forcing in connected graphs.


\section{Fundamental Properties of $F_t(G)$}
\label{S:properties}

In this section we provide fundamental properties of total forcing. Since total forcing is not defined on graphs with isolated vertices, we shall assume throughout that all graphs contain no isolates. Since any forcing set on a single vertex will induce a graph with an isolated vertex, we observe that all TF-sets must contain at least two vertices.

\begin{ob}\label{twoVertices}
If $G$ is an isolate-free graph, then $F_t(G) \ge 2$.
\end{ob}

As observed in \cite{DaHeMaPe16}, the only graphs with $F(G) = F_c(G) = 1$ are paths. Moreover, they also observed $\gamma_P(G)\le F(G)$. Indeed, if $G$ is not a path, then $F(G) \ge 2$ and $F_c(G)\ge 2$. Next observe that all TF-sets are forcing sets, and all connected forcing sets are TF-sets (provided $F_c(G)\ge 2$). We combine these observations and present them in the following chain of inequalities.

\begin{ob}\label{Forcing Relation}
If $G$ is a connected graph that is not a path, then
\[
\gamma_P(G) \le F(G) \le F_t(G) \le F_c(G).
\]
\end{ob}

We remark that Observation \ref{Forcing Relation} yields an analogous chain of inequalities on forcing to that of domination presented in Observation \ref{domOb}.

Let $S \subseteq V$ be a forcing set of an isolate-free graph $G$. For each vertex $v\in S$, color exactly one neighbor of $v$. Call this coloring $W$. This resulting set is a superset of $S$, and hence, is also a forcing set. Moreover, $G[W]$ contains no isolated vertices. Hence, $W$ is a TF-set. This observation implies that the total forcing number is no more than twice the forcing number. For example, a leaf of every non-trivial path is a forcing set of the path, while any two adjacent vertices on the path form a TF-set of the path, implying that for $n\ge 2$, $F_t(P_n) = 2$, $F(P_n) = 1$, and therefore $F_t(P_n) = 2F(P_n)$. We state our result formally as follows.

\begin{ob}\label{twiceForcing}
If $G$ is an isolate-free graph, then $F_t(G) \le 2F(G)$, and this bound is sharp.
\end{ob}

We next classify the total forcing number of simple classes of graphs. As shown in \cite{DaHeMaPe16}, $F(C_n) = F_c(C_n) =2$ and $F(K_n) = F_c(K_n) = n-1$, whenever $n\ge 3$. This together with Observation \ref{Forcing Relation} imply $F_t(C_n) = 2$ and $F_t(K_n) = n-1$ whenever $n\ge 3$. We state these results formally with the following observation.

\begin{ob}
\label{simpleFormula}
For $n\ge 3$, the following holds. \\
\indent {\rm (a)} $F_t(P_n) = 2$. \\
\indent {\rm (b)} $F_t(C_n) = 2$. \\
\indent {\rm (c)} $F_t(K_n) = n -1$.
\end{ob}

We next observe that if $G$ is a connected graph of order $n\ge 3$, then coloring any set of $n-1$ vertices while leaving an arbitrary minimum degree vertex non-colored will result in a colored subgraph without isolates. Moreover, any colored vertices under such a coloring will have at most one non-colored neighbor, and hence will be a TF-set. Note that this coloring is best possible for the complete graph $K_n$ of order $n\ge 3$. Hence, we establish the following simple upper bound.
\begin{ob}\label{simpleUpperBound}
If $G$ is an isolate-free graph of order $n\ge 3$, then $F_t(G) \le n - 1$, and this bound is sharp.
\end{ob}

In light of Observation \ref{simpleUpperBound}, it is natural to ask what graphs on a given number of vertices possess largest total forcing numbers? Before answering this question, we recall a characterization of $F_c(G) = n-1$ presented in \cite{Brimkov Davila} and \cite{DaHeMaPe16}.

\begin{thm}[\cite{Brimkov Davila,DaHeMaPe16}]\label{orderMinusOne}
If $G$ is a connected graph of order $n\ge 3$, then $F_c(G)= n -1$ with equality if and only if $G = K_n$ for $n \ge 3$, or $G = K_{1,n-1}$ with $n\ge 4$.
\end{thm}

As an immediate consequence of Observations \ref{Forcing Relation} and \ref{simpleUpperBound}, and Theorem \ref{orderMinusOne} above, we have the following result.

\begin{thm}\label{orderMinusOneIFF}
If $G$ is a connected graph of order $n\ge 3$, then $F_t(G)= n - 1$, with equality if and only if $G = K_n$ for $n \ge 3$, or $G = K_{1,n-1}$ with $n\ge 4$.
\end{thm}

We establish next a useful property of a total forcing set in graphs that contain a vertex with at least two leaf neighbors.

\begin{lem}\label{leafLem}
If $G$ is an isolate-free graph, then every vertex with at least two leaf neighbors is contained in every TF-set, and all except possibly one leaf neighbor of such a vertex is contained in every total forcing set.
\end{lem}
\proof Let $G$ be an isolate-free graph, and let $v$ be a vertex of $G$ with at least two leaf neighbors. Let $S$ be an arbitrary TF-set of $G$. If the vertex $v$ is not in $S$, then $S$ contains at least one leaf neighbor of $S$. Such a leaf neighbor of $v$ would be an isolated vertex in $G[S]$, and so the set $S$ would not be a TF-set of $G$, a contradiction. Therefore, $v\in S$. If at least two leaf neighbors of $v$ are not in $S$, then the vertex $v$ cannot be an $S$-forcing vertex, implying that no $S$-uncolored leaf neighbor of $v$ becomes colored during the forcing process, a contradiction.~\qed


\section{Total Forcing versus Domination Parameters}
\label{S:TotalFvsDom}

It is well known that both the forcing number and connected forcing number of graphs are related to various domination parameters. In this section, we explore the relationship between total forcing and various domination parameters. We begin with the following result which relates the total forcing number to the domination number.

\begin{thm}
If $G$ is a graph with minimum degree at least~$3$ and maximum degree~$\Delta$, then
\[
F_t(G) \le \gamma(G)\Delta,
\]
with equality if and only if $G$ has a perfect dominating set and every vertex in this set has degree $\Delta$ in $G$.
\end{thm}
\proof Let $G$ be a graph with minimum degree at least~$3$ and maximum degree $\Delta = \Delta(G)$, and let $S$ be a minimum dominating set of $G$, and so $|S| = \gamma(G)$. Next let $\barS = V\setminus S$, and let $W\subseteq V$ be the set of colored vertices obtained by coloring all of $S$, and for each $v\in S$, coloring all but one neighbor of $v$ in $\barS$. Observe that every vertex in $S$ has at most one $W$-uncolored neighbor in $\barS$. Further, each vertex of $S$ with exactly one $W$-uncolored neighbor in $\barS$ is a $W$-forcing vertex. Moreover, since every vertex in $G$ is either in $S$, or has a neighbor in $S$, it follows that every vertex of $G$ is either in $W$, or is forced by a vertex in $W$. Since $G$ has the property that $\delta \ge 3$, it follows that $W$ induces a graph without isolates. Hence, $W$ is a TF-set of $G$. We remark that the number of $W$-colored vertices is the number of vertices in $S$ together with at most $d_G(v) - 1$ neighbors of each vertex $v\in S$ that belong to $\barS$. Thus,
\[
\begin{array}{lcl}
F_t(G) & \le & \displaystyle{|W|} \1 \\
& \le & \displaystyle{|S| + \sum_{v \in S} (d_G(v)-1)} \1 \\
& \le & \displaystyle{|S| + |S|(\Delta -1)} \1 \\
& = & \displaystyle{|S|\Delta} \1 \\
& = & \displaystyle{\gamma(G)\Delta}.
\end{array}
\]

As shown above, $F_t(G) \le \gamma(G)\Delta$. Suppose next that $F_t(G) = \gamma(G)\Delta$. Thus we must have equality throughout the above inequality chain. In particular, $d_G(v) = \Delta$ for every vertex $v\in S$. Further, the neighbors of each vertex $v\in S$ belong to $\barS$ and we color all but one, namely $d_G(v) - 1 = \Delta -1$, such neighbors. If two vertices, say $u$ and $v$, in $S$ have a common neighbor, $x$ say, then we color $\Delta - 1$ neighbors of $u$ including the vertex $x$, and we color $\Delta -2$ neighbors of $v$ excluding the vertex $x$ and one other arbitrary neighbor of $v$. Coloring $\Delta -1$ neighbors of all other vertices in $S$ different from $u$ and $v$ produces a TF-set of $G$ with cardinality strictly less than $|S| + |S|(\Delta -1)$, noting that the forcing vertex $u$ is played in the time step before the forcing vertex $v$ is played. This produces a TF-set of size less than $\gamma(G)\Delta$, a contradiction. Therefore, no two vertices in $S$ are adjacent, or have a common neighbor. Thus, $S$ is a perfect dominating set, and every vertex in $S$ has degree $\Delta$ in $G$.~\qed

\medskip
We next relate the total forcing number to the total domination number.

\begin{thm}
If $G$ is an isolate-free graph with maximum degree $\Delta \ge 2$, then
\[
F_t(G) \le \gamma_t(G)(\Delta -1),
\]
with equality if and only if $G$ has a perfect total dominating set and every vertex in this set has degree $\Delta$ in $G$.
\end{thm}
\proof Let $G$ be an isolate-free graph with maximum degree $\Delta = \Delta(G) \ge 2$. Let $S\subseteq V$ be a minimum total dominating set of $G$, and so $|S| = \gamma_t(G)$. Let $\barS = V \setminus S$. Next let $W$ be a set of colored vertices obtained by coloring $S$, and for each $v\in S$, coloring all but one neighbor of $v$ in $\barS$. Observe that every vertex in $S$ has at most one $W$-uncolored neighbor in $\barS$. Further, each vertex of $S$ with exactly one $W$-uncolored neighbor in $\barS$ is a $W$-forcing vertex. Moreover, since every vertex in $G$ is either in $S$, or has a neighbor in $S$, it follows that every vertex of $G$ is either in $W$, or is forced by a vertex in $W$. Since each vertex of $S$ has at least one neighbor in $S$, and we have only colored vertices in $\barS$ adjacent with vertices in $S$, we observe that $W$ is a TF-set. We remark that the number of $W$-colored vertices is the number of vertices in $S$ together with at most $d_G(v)-2$ neighbors of each vertex $v\in S$ that belong to $\barS$. Thus,
\[
\begin{array}{lcl}
F_t(G) & \le & \displaystyle{|W|} \1 \\
& \le & \displaystyle{|S| + \sum_{v \in S} (d_G(v)-2)} \1 \\
& \le & \displaystyle{|S| + |S|(\Delta -2)} \1 \\
& = & \displaystyle{|S|(\Delta-1)} \1 \\
& = & \displaystyle{\gamma_t(G)(\Delta -1)}.
\end{array}
\]

As shown above, $F_t(G)\le \gamma_t(G)(\Delta - 1)$. Next suppose that $F_t(G) = \gamma_t(G)(\Delta - 1)$. Thus, we must have equality throughout the above inequality chain. In particular, $d_G(v) = \Delta$ for every $v\in S$. Further, each vertex $v\in S$ has all, except for exactly one, of its neighbors in $\barS$ and we color all but one, namely $d_G(v) - 2 = \Delta -2$, such neighbors. If two vertices, say $u$ and $v$, in $S$ have a common neighbor, $x$ say, in $\barS$, then we color $\Delta -2$ neighbors of $u$ in $\barS$ including the vertex $x$, and we color $\Delta - 3$ neighbors of $v$ in $\barS$ excluding the vertex $x$ and one other arbitrary neighbor of $v$ in $\barS$. Coloring $\Delta - 2$ neighbors of all other vertices in $S$ (different from $u$ and $v$) that belong to $\barS$ produces a TF-set of $G$ with cardinality strictly less that $|S| + |S|(\Delta -2)$, noting that the forcing vertex $u$ is played in the time step before the forcing vertex $v$ is played. This produces a TF-set of cardinality less than $\gamma_t(G)(\Delta -1)$, a contradiction. Therefore, no two vertices in $S$ have a common neighbor. Thus, the set $S$ is a perfect total dominating set and every vertex in $S$ has degree $\Delta$ in $G$.~\qed

\medskip
The following result shows that the sum of the total forcing number and the (distance) $2$-domination number of a graph is at most its order.

\begin{thm}
If $G$ is a graph of order $n$ with minimum degree at least~$2$, then
\[
F_t(G) + \gamma_2(G) \le n.
\]
\end{thm}
\proof Let $G$ be a graph of order $n$ with $\delta(G) \ge 2$. Let $S$ be a minimum TF-set, and so $|S| = F_t(G)$. We show that the set $V(G) \setminus S$ is a $2$-dominating set of $G$. If this is not the case, then there is a vertex $v \in S$ at distance at least~$3$ from every vertex outside $S$. Thus, every vertex within distance~$2$ from $v$ belongs to the set $S$. We now consider the set $S' = S \setminus \{v\}$. If $G[S']$ contains an isolated vertex $w$, then since $G[S]$ is isolate-free this would imply that $v$ is the only neighbor of $w$ that belongs to the set $S$. However, the minimum degree at least two condition implies that $w$ has at least one neighbor outside $S$. Such a neighbor of $w$ that belongs to $V(G) \setminus S$ is at distance~$2$ from $v$ in $G$, a contradiction. Therefore, $G[S']$ is isolate-free. Further since $S$ is a TF-set and since the vertex $v$ is never played in the forcing process, the set $S'$ is a forcing set. Hence, the set $S'$ is a TF-set, contradicting the minimality of the set $S$. We deduce, therefore, that the set $V(G) \setminus S$ is a $2$-dominating set of $G$. Thus, $\gamma_2(G) \le |V(G) \setminus S| = n - |S| = n - F_t(G)$.~\qed

\medskip
We next recall a lemma which relates forcing sets to power dominating sets.

\begin{lem}[\cite{Davila,DaHeMaPe16}]
\label{powerDomLem}
Let $G$ be a graph. Then, $S\subseteq V$ is a power dominating set if and only if $N[S]$ is a forcing set of $G$.
\end{lem}

Using Lemma \ref{powerDomLem}, we next establish two relationships between the total forcing number and the power domination number.

\begin{thm}
\label{p:power}
If $G$ is an isolate-free graph with maximum degree $\Delta$, then
\[
F_t(G) \le \gamma_P(G)(\Delta +1).
\]
\end{thm}
\proof Let $G$ be an isolate-free graph with maximum degree $\Delta = \Delta(G)$. Let $S\subseteq V$ be a minimum power dominating set of $G$; that is, $|S| = \gamma_P(G)$. By the definition of power domination, see property $(1)$ in the definition, we observe that $S$ power dominates (observes) its closed neighborhood. Moreover, since $G$ contains no isolated vertices, we observe that $\cP_G^0 = N_G[S]$ is a set of vertices which induces a graph without isolated vertices. By Lemma \ref{powerDomLem}, it follows that $\cP_G^0 = N_G[S]$ is a forcing set. In particular, if we color $W =\cP_G^0 = N_G[S]$, we obtain a TF-set. We remark the number of $W$-colored vertices is precisely the number of vertices in $N[S]$, i.e., we color all the vertices of $S$ and at most $d_G(v)$ neighbors of each $v\in S$ that belong to $\barS$. Thus,

\[
\begin{array}{lcl}
F_t(G) & \le & \displaystyle{|W|} \1 \\
& \le & \displaystyle{|S| + \sum_{v \in S}d_G(v)} \1 \\
& \le & \displaystyle{|S| + |S|\Delta} \1 \\
& = & \displaystyle{|S|(\Delta + 1)} \1 \\
& = & \displaystyle{\gamma_P(G)(\Delta +1)}. \hspace*{1cm} \Box
\end{array}
\]

\medskip
The following fundamental result in domination is attributed to Ore \cite{Ore's Theorem}, and is often referred to as Ore's Theorem.

\begin{thm}[Ore's Theorem \cite{Ore's Theorem}]\label{OresThm}
If $G$ is an isolate-free graph of order $n$, then $\gamma(G)\le \frac{n}{2}$.
\end{thm}

Making use of Ore's Theorem, we next present another relationship between the power domination number and the total forcing number.

\begin{thm}\label{powerDom}
If $G$ is an isolate-free graph, then $F_t(G) \ge 2\gamma_P(G)$, and this bound is sharp.
\end{thm}
\proof Let $G$ be an isolate-free graph and let $S \subseteq V$ be a $F_t(G)$-set. Since $S$ is a TF-set, we observe that $G[S]$ is an isolate-free graph. Thus, by Theorem \ref{OresThm}, there exists a set of vertices $W \subseteq S$ which dominates $G[S]$ and has cardinality at most $|S|/2$. Since $W$ dominates $G[S]$, it follows that $S\subseteq N[W]$. By Lemma \ref{powerDomLem}, $W$ is a power dominating set of $G$, and so $\gamma_P(G) \le |W| \le \frac{1}{2}|S| = \frac{1}{2}F_t(G)$. Rearranging, we get our desired result. To see that this bound is sharp, consider the cycle $C_n$ on $n$ vertices.~\qed


\section{Effects of Vertex Removal on Total Forcing}
\label{S:vertexRem}

If $v$ is a vertex of an isolate-free graph $G$ whose removal creates no isolates, it is natural to consider the effect on the total forcing number when the vertex $v$ is deleted from $G$. In this section, we discuss how the total forcing number of such a graph $G$ relates to the total forcing number of $G-v$. We answer this with the following result.
\begin{thm}
\label{t:vertex}
If $G$ is an isolate-free graph with maximum degree $\Delta$ that contains a vertex $v$ such that $G-v$ contains no isolates, then
\begin{equation*}
F_t(G) - 2 \le F_t(G-v) \le F_t(G) + \Delta,
\end{equation*}
and these bounds are tight.
\end{thm}
\proof Let $v$ be a vertex in the isolate-free graph of maximum degree $\Delta$ such that $G-v$ contains no isolates. We note that $\Delta \ge 2$. Let $T \subseteq V(G) \setminus \{v\}$ be a minimum TF-set of $G-v$, and so $|T| = F_t(G-v)$. Further, let $w$ be an arbitrary neighbor of $v$ in $G$. The set $T \cup\{v,w\}$ is a TF-set of $G$, implying that $F_t(G) \le |T \cup\{v,w\}| \le |T| + 2 = F_t(G-v) + 2$. Rearranging, we obtain $F_t(G) - 2 \le F_t(G-v)$. This establishes the desired lower bound on $F_t(G-v)$.

We next prove the upper bound. Let $S \subseteq V(G)$ be a minimum TF-set of $G$, and so $|S| = F_t(G)$. Since $S$ is a TF-set of $G$, the graph $G[S]$ contains no isolated vertex. Suppose that $v \notin S$. Immediately before the vertex $v$ is played in the forcing process, at most one neighbor of $v$ is uncolored. If immediately before the vertex $v$ is played, no neighbor of $v$ is uncolored, then the set $S \setminus \{v\}$ is TF-set of $G - v$, implying that $F_t(G-v) \le |S| - 1 = F_t(G) - 1$. If immediately before the vertex $v$ is played, exactly one neighbor, say $w$, of $v$ is uncolored, then let $w'$ denote an arbitrary neighbor of $w$ different from $v$. In this case, the set $(S \setminus \{v\}) \cup \{w,w'\}$ is a TF-set of $G - v$, implying that $F_t(G-v) \le (|S| - 1) + 2 = F_t(G) + 1 \le F_t(G) + \Delta - 1$.

Hence, we may assume that $v \in S$, for otherwise the desired upper bound follows. Let $S_v$ denote the set of neighbors of $v$ in $S$, and so $S_v = N(v) \cap S$. By assumption, each neighbor of $v$ has degree at least~$2$ in $G$. For each vertex $u \in S_v$, let $u'$ be a neighbor of $u$ different from $v$. We now consider the set,
\[
S_v' = \bigcup_{u \in S_v} \{u'\}.
\]

If $N_G(v) = S_v$ or if the vertex $v$ is not played in the forcing process, then since $S$ is a TF-set of $G$, the set $(S \setminus \{v\}) \cup S_v'$ is a TF-set of $G-v$, implying that $F_t(G-v) \le |S| + |S_v'| - 1 \le |S| + \Delta - 1 = F_t(G) + \Delta - 1 < F_t(G) + \Delta$. If the vertex $v$ is played in the forcing process, then let $v_1$ denote the uncolored neighbor of $v$ that becomes colored when $v$ is played and let $v_2$ be an arbitrary neighbor of $v_1$ different from $v$. In this case, we note that at least one neighbor of $v$, namely $v_1$, does not belong to the set $S$, and so $|S_v'| \le d_G(v) - 1 \le \Delta - 1$. Further, the set $(S \setminus \{v\}) \cup S_v' \cup \{v_1,v_2\}$ is a TF-set of $G-v$, implying that $F_t(G-v) \le (|S| - 1) + |S_v'| + 2  \le |S| + \Delta = F_t(G) + \Delta$.

That the upper bound of Theorem~\ref{t:vertex} is tight may be seen as follows. For $k \ge 3$, if $G$ is the graph obtained by subdividing each edge of the star $K_{1,k}$ exactly once, then $F_t(G) = k$ and $\Delta(G) = k$. However, if $v$ is the unique (central) vertex of degree~$k$ in $G$, then $G-v$ is the graph consisting of $k$ disjoint copies of $P_2$. Thus, $F_t(G-v) = 2k = F_t(G) + \Delta$. When $k = 3$, the graph $G$ is illustrated in Figure~\ref{f:fig1}(a).

That the lower bound of Theorem~\ref{t:vertex} is tight may be seen by taking, for example, $G$ to be the graph illustrated in Figure~\ref{f:fig1}(b) where the vertex $v$ is an arbitrary leaf. In this example, $F_t(G) = 5$ and $F_t(G-v) = 3 = F_t(G) - 2$. We remark that there are connected graphs $G$ of arbitrarily large order achieving the lower bound of Theorem~\ref{t:vertex}. For example, taking $k \ge 2$ vertex disjoint copies $G_1, G_2, \ldots, G_k$ of the graph $G$ illustrated in Figure~\ref{f:fig1}(b), where $w_i$ denotes the vertex of $G_i$ named $w$ in $G$ for $i \in [k]$. Let $H$ be the graph formed by the disjoint union of these $k$ graphs by adding the edges $w_jw_{j+1}$ for $j \in [k-1]$. We note that $F_t(H) = 5k$ and $F_t(H - x) = 5k - 2 = F_t(H) - 2$ for every leaf $x$ in $H$.~\qed

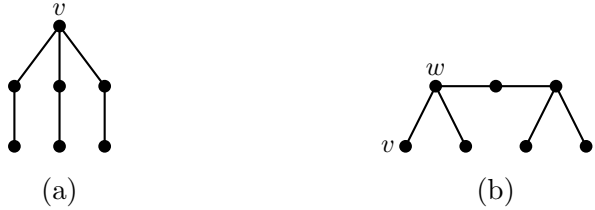
\begin{figure}[htb]
\begin{center}
\begin{tikzpicture}[scale=.8,style=thick,x=1cm,y=1cm]
\def\vr{2.5pt} 
%
\path (0,0) coordinate (w1);
\path (0.75,0) coordinate (w2);
\path (1.5,0) coordinate (w3);
\path (0,1) coordinate (u1);
\path (0.75,1) coordinate (u2);
\path (1.5,1) coordinate (u3);
\path (0.75,2) coordinate (v);
%
\draw (v) -- (u1);
\draw (v) -- (u2);
\draw (v) -- (u3);
\draw (u1) -- (w1);
\draw (u2) -- (w2);
\draw (u3) -- (w3);
\draw (v) [fill=black] circle (\vr);
\draw (u1) [fill=black] circle (\vr);
\draw (u2) [fill=black] circle (\vr);
\draw (u3) [fill=black] circle (\vr);
\draw (w1) [fill=black] circle (\vr);
\draw (w2) [fill=black] circle (\vr);
\draw (w3) [fill=black] circle (\vr);
\draw (0.75,-0.75) node {(a)};
\draw[anchor = south] (v) node {{\small $v$}};
\path (6.5,0) coordinate (u1);
\path (7.5,0) coordinate (u2);
\path (8.5,0) coordinate (u5);
\path (9.5,0) coordinate (u6);

\path (7,1) coordinate (v1);
\path (8,1) coordinate (v2);
\path (9,1) coordinate (v3);

\draw (v1) -- (u1);
\draw (v1) -- (u2);
\draw (v3) -- (u5);
\draw (v3) -- (u6);
\draw (v1) -- (v3);

\draw (v1) [fill=black] circle (\vr);
\draw (v2) [fill=black] circle (\vr);
\draw (v3) [fill=black] circle (\vr);
\draw (u1) [fill=black] circle (\vr);
\draw (u2) [fill=black] circle (\vr);
\draw (u5) [fill=black] circle (\vr);
\draw (u6) [fill=black] circle (\vr);
\draw[anchor = east] (u1) node {{\small $v$}};
\draw[anchor = south] (v1) node {{\small $w$}};
\draw (8,-0.75) node {(b)};
\end{tikzpicture}
\end{center}
\vskip -0.5 cm
\caption{Graphs illustrating tightness of Theorem~\ref{t:vertex}} \label{f:fig1}
\end{figure}

We next show next that the upper bound of Theorem~\ref{t:vertex} can be improved if we remove a vertex of degree~$1$. More precisely, we show that removing a vertex of degree~$1$ from a connected graph of order at least~$3$ cannot increase the total forcing number.

\begin{lem}\label{treeSpread}
If $G$ is a connected graph of order at least~$3$ and $v$ a vertex of degree~$1$ in~$G$, then \[F_t(G-v) \le F_t(G).\]
\end{lem}
\proof Let $G$ is a connected graph of order at least~$3$. Let $v$ be a vertex of degree~$1$ in $G$ and let $w$ be the neighbor of $v$. Since $G$ has order at least~$3$, we note that $w$ has at least two neighbors and that $G-v$ is a non-trivial graph. Let $S$ be a minimum TF-set of $G$, and so $|S| = F_t(G)$. If $v \notin S$, then the set $S$ is also a TF-set of $G - v$, and so $F_t(G-v) \le |S| = F_t(G)$. Hence, we may assume that $v \in S$, for otherwise the desired result follows. Since $G[S]$ contains no isolated vertex, $w \in S$. If every neighbor of $w$ belongs to $S$, then $S \setminus \{v\}$ is a TF-set of $G$, contradicting the minimality of $S$. Hence, there is a neighbor $x$ of $w$ that does not belong to $S$. Replacing $v$ in $S$ with the vertex $x$ produces a new minimum TF-set, $S'$ say, of $G$ that is also a TF-set of $G - v$. Thus, $F_t(G-v) \le |S'| = |S| = F_t(G)$.~\qed


\section{Complexity of Total Forcing}
\label{S:NP}

As shown in the previous section, the total forcing number is related to a myriad of $NP$-complete graph invariants. In this section, we show that the decision problem associated with total forcing is also $NP$-complete. We state this decision problem formally as follows.

\noindent
PROBLEM: Total forcing (TF)
\\
INSTANCE: An isolate-free graph $G =(V,E)$ of order $n$ and a positive integer $k\le n$.
\\QUESTION: Does there exist a total forcing set $S\subseteq V$ of cardinality at most~$k$?

\medskip
With the following theorem we prove that total forcing is $NP$-complete.

\begin{thm}\label{NP}
TF is $NP$-complete
\end{thm}
\proof We first establish that TF is $NP$. Let $G$ be an isolate-free graph of order $n$, and let $S\subseteq V$. Next observe that it can be checked in polynomial time if a vertex $v\in S$ has exactly one neighbor which is not in $S$; that is, we may check if $v$ is a $S$-forcing vertex in polynomial time. Next observe that there are at most $n-1$ forcing steps during the forcing process in $G$. Hence, there exists a nondeterministic algorithm which may check in polynomial time whether $S$ is a forcing set, whether $S$ induces a graph without isolated vertices, and if $S$ has cardinality at most $k$. It follows that TF is $NP$.

Next, we use for our reduction, the problem of forcing, which has been shown to be $NP$-complete in \cite{zf_np}. We state the decision problem of forcing below.

\noindent
PROBLEM: Forcing ($ZF$)
\\INSTANCE: A simple graph $G =(V,E)$ of order $n$ and a positive integer $k\le n$.
\\QUESTION: Does there exist a forcing set of vertices $S\subseteq V$ with cardinality at most~$k$?

In order to prove that TF is $NP$-complete we construct a transformation from $ZF$ to TF. Let $I = \langle G, k \rangle$ be an instance of $ZF$, where $G = (V,E)$ and $V = \{v_1, \dots, v_n\}$. Let $G'$ be the graph obtained from $G$ as follows. For each vertex $v$ of $G$, add a vertex disjoint copy of a path $P_3$ and an edge from $v$ to the (central) vertex of degree~$2$ in $P_3$. For each vertex $v_i$ of $G$, let $G_i$ be the added copy of $P_3$ associated with $v_i$ for $i\in [n]$. Further, let $G_i$ be given the path $l_{i}^1v_{i}^*l_{i}^2$, and so $v_{i}v_{i}^{*}$ is an edge of $G'$. The construction of the graph $G'$ from the graph $G$ is illustrated in Figure~\ref{f:fig1}. We next define $f(I) = \langle G', k + 2n \rangle$.


\begin{figure}[htb]
\begin{center}
\begin{tikzpicture}[scale=.8,style=thick,x=1cm,y=1cm]
\def\vr{2.5pt} 
\path (0,0.5) coordinate (x1);
\path (1.5,0.5) coordinate (x2);
\path (3,0.5) coordinate (x3);
\path (4.5,0.5) coordinate (x4);
\path (6,0.5) coordinate (x5);
\path (0,1.25) coordinate (y1);
\path (1.5,1.25) coordinate (y2);
\path (3,1.25) coordinate (y3);
\path (4.5,1.25) coordinate (y4);
\path (6,1.25) coordinate (y5);
\path (-0.4,2) coordinate (w1);
\path (0.4,2) coordinate (z1);
\path (1.1,2) coordinate (w2);
\path (1.9,2) coordinate (z2);
\path (2.6,2) coordinate (w3);
\path (3.4,2) coordinate (z3);
\path (4.1,2) coordinate (w4);
\path (4.9,2) coordinate (z4);
\path (5.6,2) coordinate (w5);
\path (6.4,2) coordinate (z5);
\draw (2.25,0.1) node {$G$};
\draw (-1.5,1) node {$G'$};
%
\draw (x1) -- (y1);
\draw (x2) -- (y2);
\draw (x3) -- (y3);
\draw (x4) -- (y4);
\draw (x5) -- (y5);
\draw (y1) -- (w1);
\draw (y1) -- (z1);
\draw (y2) -- (w2);
\draw (y2) -- (z2);
\draw (y3) -- (w3);
\draw (y3) -- (z3);
\draw (y4) -- (w4);
\draw (y4) -- (z4);
\draw (y5) -- (w5);
\draw (y5) -- (z5);
\draw (x1) [fill=black] circle (\vr);
\draw (x2) [fill=black] circle (\vr);
\draw (x3) [fill=black] circle (\vr);
\draw (x4) [fill=black] circle (\vr);
\draw (x5) [fill=black] circle (\vr);
\draw (y1) [fill=black] circle (\vr);
\draw (y2) [fill=black] circle (\vr);
\draw (y3) [fill=black] circle (\vr);
\draw (y4) [fill=black] circle (\vr);
\draw (y5) [fill=black] circle (\vr);
\draw (w1) [fill=black] circle (\vr);
\draw (z1) [fill=black] circle (\vr);
\draw (w2) [fill=black] circle (\vr);
\draw (z2) [fill=black] circle (\vr);
\draw (w3) [fill=black] circle (\vr);
\draw (z3) [fill=black] circle (\vr);
\draw (w4) [fill=black] circle (\vr);
\draw (z4) [fill=black] circle (\vr);
\draw (w5) [fill=black] circle (\vr);
\draw (z5) [fill=black] circle (\vr);
\draw [rounded corners] (-0.5,-0.25)
rectangle (6.65,0.7) node [black,right] {};
\end{tikzpicture}
\end{center}
\vskip -0.25cm
\caption{Obtaining $G'$ from $G$} \label{f:fig1}
\end{figure}
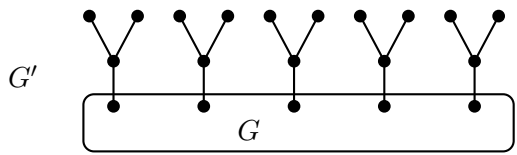

Next suppose that $I =\langle G, k \rangle$ is a ``yes" instance of $ZF$. That is, $G$ has a forcing set $S\subseteq V$ of cardinality at most $k$. We next show that $S' = S \cup\{v_1^*,\dots, v_n^*, l_1, \dots, l_n\}$ is a TF-set of $G'$. It is clear that $G'$ is isolate-free since if $v$ was isolated in $G$, it is now adjacent to the degree two vertex of some copy of $P_3$. Moreover, each vertex $S$ has a neighbor which is colored in some copy of $P_3$, and hence, $S'$ induces a graph with no isolated vertices. By Lemma~\ref{leafLem}, it is also clear that each $v_i^*$, and all but one of their respective degree one neighbors, at least $l_i$ (say), are contained in every TF-set of $G'$, for $1\le i \le n$; of which $S'$ satisfies. Next observe that since each $v_i^*$ is colored, any colored vertex in $G$, at any point during the forcing process starting in $S$, will have the same number of non-colored neighbors in $G'$. That is, we are assured that $G$ as a subgraph will be completely vertex colored as a subgraph due to the forcing process starting at $S'$ in $G'$. Moreover, once the vertices of $G$ are colored in $G'$, each $v_i^*$, for $1\le i \le n$, will have exactly one non-colored neighbor which is a leaf different that $l_i$, and hence will force. It follows that $S'$ is a TF-set in $G'$. It follows that $S'$ is a TF-set of $G'$ with cardinality at most $k+2n$, and hence, $f(I) = \langle G', k + 2n\rangle$ is a ``yes" instance of TF.

Conversely, suppose $f(I) = \langle G', k + 2n \rangle$ is a ``yes" instance of TF. That is, $G'$ has a TF-set $S'$ of cardinality $k + 2n$. By Lemma \ref{leafLem}, we may assume $\{v_1^*, \dots, v_n^*, l_1, \dots, l_n\} \subseteq S'$. With this observation, we next observe for $1\le i \le n$, it is the case that $v_i^*$ is the unique vertex which forces its degree one neighbor different from $l_i$ in $G'$. In particular, during the forcing process in $G'$ any vertex which is colored in the subgraph $G$ is colored by a vertex contained in $G$, i.e., the set $\{v_1^*, \dots, v_n^*, l_1, \dots, l_n\} \subseteq S'$ forces no vertex of $G$. It follows that $S=S'\setminus \{v_1^*, \dots, v_n^*, l_1, \dots, l_n\}$ is a forcing set of $G$. Since $S'$ has cardinality at most $k + 2n$, it follows that $S$ has cardinality at most $k$. In particular, we have that $I$ is a ``yes" instance of ZF. \qed

\section{A General Upper Bound}
\label{S:bds}

In this section we provide a sharp upper bound on the total forcing number of a graph in terms of its maximum degree and order. Recall that for a subset $S \subseteq V(G)$, the degree of $v$ in $S$, denoted $d_S(v)$, is the number of vertices in $S$ adjacent to~$v$. We are now in a position to prove the following upper bound on the total forcing number of a graph.

\begin{thm}
\label{t:upperbd}
If $G$ is a connected graph of order $n \ge 3$ with maximum degree~$\Delta$, then
\[
F_t(G) \le \left( \frac{\Delta}{\Delta + 1} \right) n,
\]
with equality if and only if $G \cong K_{\Delta + 1}$ or $G \cong K_{1,\Delta}$.
\end{thm}
\proof Let $G$ be a connected graph of order $n \ge 3$  with maximum degree~$\Delta$. If $\Delta = 2$, then $G \cong P_n$ or $G \cong C_n$. In both cases, $F_t(G) = 2 \le \frac{2}{3} n = ( \frac{\Delta}{\Delta + 1} ) n$, as desired. Further, if $F_t(G) = ( \frac{\Delta}{\Delta + 1} ) n$, then we must have equality throughout this inequality chain, implying that $n = 3$ and $G \cong C_3 = K_3$ or $G \cong P_3$.
Hence, we may assume that $\Delta \ge 3$. Among all maximum packings in $G$, let $P = \{v_1, \ldots, v_k\}$ be chosen so that

\indent (1) the number of vertices not dominated by $P$ is a maximum,
\\
\indent (2) subject to (1), the sum of degrees of vertices in $P$ is a maximum.

We note that if $P$ is a perfect packing, then every vertex of $G$ is dominated by $P$. Let $A$ be the set of all neighbors of vertices in $P$; that is, $A = N(P)$. For $i \in [k]$, let $A_i = N(v_i)$, and so $A = \cup_{i=1}^k A_i$. Further, let $|A_i| = \alpha_i$, and so
$1 \le \alpha_i = d_G(v_i) \le \Delta$. Since $P$ is a packing, we note that $(A_1,\ldots,A_k)$ is a partition of $A$. Further, $(A_1 \cup \{v_1\},\ldots,A_k \cup \{v_k\})$ is a partition of $N[P]$.

Let $B = V(G) \setminus N[P]$. If a vertex $w \in B$ is at distance at least~$3$ from every vertex of $P$, then the set $P \cup \{w\}$ would be a packing in $G$, contradicting the maximality of $P$. Hence, every vertex in $B$ has a neighbor in $A$ and is therefore adjacent to a vertex in $A_i$ for some $i \in [k]$. Equivalently, the set $A$ dominates the set $B$.

We now define a weak partition $(B_1, \ldots, B_k)$ of the set $B$ as follows, where by a weak partition we mean a partition where some of the sets may be empty. Let $B_1$ be the set of all vertices in $B$ that have a neighbor in $A_1$. Thus, $B_1$ consists of all vertices $v \in B$ such that $d_{A_1}(v) \ge 1$. For $i \in [k] \setminus \{1\}$, let $B_i$ be the set of all vertices in $B$ that have a neighbor in $A_i$ but no neighbor in $A_j$ for any $j \in [i-1]$. Thus for $i \in [k] \setminus \{1\}$, the set $B_i$ consists of all vertices $v \in B$ such that $d_{A_j}(v) = 0$ for all $j$ where $1 \le j < i$ and $d_{A_i}(v) \ge 1$. By definition, $(B_1, \ldots, B_k)$ is a weak partition of the set $B$.

If $d_G(v_i) = 1$ and $B_i = \emptyset$ for all $i \in [k]$, then $V(G) = P \cup A$, $|A| = |P| = k$ and $n = 2k$. In this case, by the connectivity of $G$ and since each vertex of $A$ has exactly one neighbor outside $A$, the set $A$ is a TF-set of $G$, implying that
\[
F_t(G) \le |A| = k = \frac{1}{2}n < \frac{3}{4}n \le \left( \frac{\Delta}{\Delta + 1} \right) n.
\]

Hence, we may assume that for at least one $i \in [k]$, $d_G(v_i) \ge 2$ or $d_G(v_i) = 1$ and $B_i \ne \emptyset$. Renaming the vertices in $P$ if necessary, we may further assume that $d_G(v_i) \ge 2$ or $d_G(v_i) = 1$ and $B_i \ne \emptyset$ for all $i \in [r]$ where $1 \le r \le k$. Possibly, $r = k$. If $r < k$, then we note that $d_G(v_j) = 1$ and $B_j = \emptyset$ for all $j$ where $r+1 \le j \le k$.

For $i \in [k]$, let $G_i$ be the graph induced by $N[v_i] \cup B_i$, where recall that $N[v_i] = \{v_i\} \cup A_i$. Further, let $G_i$ have order~$n_i$. We first restrict our attention to those graphs $G_i$ where $i \in [r]$. By definition of the sets $A_i$ and $B_i$, we note that the set $A_i$ dominates the set $B_i$. Let $D_i$ be a minimum set of vertices in $A_i$ that dominate $B_i$. We now consider two cases. In both cases, we note that necessarily $i \in [r]$.

\newpage
\emph{Case~1. $B_i \ne \emptyset$.} Let $|D_i| = d_i$, and note that $1 \le d_i \le \alpha_i$. Further, we note that if $D_i = A_i$, then $d_i = \alpha_i$, while if $D_i \subset A_i$, then $d_i < \alpha_i$.
By the minimality of the set $D_i$, each vertex in $D_i$ dominates a vertex in $B_i$ that is not dominated by the other vertices in $D_i$. For each vertex $x \in D_i$, let $x'$ be a vertex in $B_i$ that is dominated in $G_i$ by $x$ but by no vertex of $D_i \setminus \{x\}$. Further, let
\[
D_i' = \bigcup_{x \in D_i} \{x'\} \hspace*{0.5cm} \mbox{and} \hspace*{0.5cm} L_i = B_i \setminus D_i'.
\]

Let $|L_i| = \ell_i$, and so $|B_i| = d_i + \ell_i$ and
\[
n_i = |A_i| + |B_i| + 1 = \alpha_i + d_i + \ell_i + 1.
\]
Each vertex in $D_i$ is adjacent to $v_i$ and to exactly one vertex in $D_i'$, and is therefore adjacent to at most~$\Delta - 2$ vertices in $L_i$, implying that
\[
\ell_i \le |D_i|(\Delta - 2) = d_i(\Delta - 2).
\]
We now consider two subcases.

\medskip
\emph{Case~1.1. $D_i \subset A_i$.} Let $|D_i| = d_i$. In this case, $1 \le d_i < \alpha_i$. In particular, we note that $\alpha_i \ge 2$. Let $w_i$ be an arbitrary neighbor of $v_i$ that does not belong to $D_i$, and let $S_i = V(G_i) \setminus (D_i' \cup \{w_i\})$. By construction, the graph $G_i[S_i]$ is isolate-free. Further, the set $S_i$ is a forcing set of $G_i$ since if $x_1 = v_i$ and $D_i = \{x_2, \ldots, x_{d_i + 1}\}$, then the sequence $x_1,x_2,\ldots,x_{d_i + 1}$ of played vertices in the forcing process results in all vertices of $G_i$ colored, where $x_i$ denotes the forcing vertex played in the $i$th step of the process. In particular, we note that $x_1 = v_i$ is an $S_i$-forcing vertex that forces $w_i$ to be colored, while each vertex in $D_i$ is an $S_i$-forcing vertex that forces its unique neighbor in $D_i'$ to be colored. Thus, the set $S_i$ is a TF-set of $G_i$. Further, $|S_i| = \alpha_i + \ell_i$. As observed earlier, $0 < d_i < \alpha_i \le \Delta$ and $\ell_i \le d_i(\Delta - 2)$, implying that $\alpha_i + \ell_i \le d_i \Delta + \alpha_i - 2d_i < d_i \Delta + \alpha_i \le \Delta(d_i + 1)$. We note that the following holds.

\[
  \begin{array}{lrcl}
  &  |S_i| & < & \left( \frac{\Delta}{\Delta + 1} \right) n_i \\
\Updownarrow & & & \\
  &  \alpha_i + \ell_i & < & \left( \frac{\Delta}{\Delta + 1} \right) ( \alpha_i + d_i + \ell_i + 1 ) \\
\Updownarrow & & & \\
  &  (\Delta + 1) (\alpha_i + \ell_i)  & < & \Delta ( \alpha_i + d_i + \ell_i + 1 ) \\
\Updownarrow & & & \\
  &  \alpha_i + \ell_i   & < & \Delta (d_i + 1). \\
  \end{array}
\]

Thus since $\alpha_i + \ell_i < \Delta (d_i + 1)$ holds, so too does the inequality $|S_i| < ( \frac{\Delta}{\Delta + 1} ) n_i$ hold. Thus, $S_i$ is a TF-set of $G_i$ satisfying the inequality
\begin{equation}
\label{Eq1}
|S_i| < \left( \frac{\Delta}{\Delta + 1} \right) n_i.
\end{equation}

\emph{Case~1.2. $D_i = A_i$.} In this case, $d_i = \alpha_i$, and so
$n_i = 2\alpha_i + \ell_i + 1$ and $\ell_i \le \alpha_i(\Delta - 2)$. Let $S_i = V(G_i) \setminus D_i'$. By construction, the graph $G_i[S_i]$ is isolate-free. Further, the set $S_i$ is a forcing set of $G_i$ since if $D_i = \{x_1, \ldots, x_{\alpha_i}\}$, then the sequence $x_1,\ldots,x_{\alpha_i}$ of played vertices in the forcing process results in all vertices of $G_i$ colored, where $x_i$ denotes the forcing vertex played in the $i$th step of the process. In particular, we note that each vertex in $D_i$ is an $S_i$-forcing vertex that forces its unique neighbor in $D_i'$ to be colored. Thus, the set $S_i$ is a TF-set of $G_i$. Further, $|S_i| = \alpha_i + \ell_i + 1$. As observed earlier, $\alpha_i \ge 1$ and $\ell_i \le \alpha_i(\Delta - 2)$, implying that $\alpha_i + \ell_i + 1 \le \Delta \alpha_i - \alpha_i + 1 \le \Delta \alpha_i$. We note that the following holds.

\[
  \begin{array}{lrcl}
  &  |S_i| & \le & \left( \frac{\Delta}{\Delta + 1} \right) n_i \\
\Updownarrow & & & \\
  &  \alpha_i + \ell_i + 1 & \le & \left( \frac{\Delta}{\Delta + 1} \right) ( 2\alpha_i + \ell_i + 1 ) \\
\Updownarrow & & & \\
  &  (\Delta + 1) (\alpha_i + \ell_i + 1)  & \le & \Delta ( 2\alpha_i + \ell_i + 1 ) \\
\Updownarrow & & & \\
  &  \alpha_i + \ell_i + 1  & \le & \Delta \alpha_i. \\
  \end{array}
\]

Thus since $\alpha_i + \ell_i + 1 \le \Delta \alpha_i$ holds, so too does the inequality $|S_i| \le ( \frac{\Delta}{\Delta + 1} ) n_i$ hold. Thus, $S_i$ is a TF-set of $G_i$ satisfying the inequality
\begin{equation}
\label{Eq2}
|S_i| \le \left( \frac{\Delta}{\Delta + 1} \right) n_i.
\end{equation}

\medskip
\emph{Case~2. $B_i = \emptyset$ and $\alpha_i \ge 2$.} We note that in this case $i \in [r]$. We now let $w_i$ be an arbitrary neighbor of $v_i$ and let $S_i = N[v_i] \setminus \{w_i\}$. Since $\alpha_i \ge 2$, we note that $G_i[S_i]$ is isolate-free. Further since $B_i = \emptyset$, we have $n_i = \alpha_i + 1$. The resulting set $S_i$ is a TF-set in $G_i$ of size
\[
|S_i| = \alpha_i = \left( \frac{\alpha_i}{\alpha_i + 1} \right) (\alpha_i + 1) \le \left( \frac{\Delta}{\Delta + 1} \right) (\alpha_i + 1) = \left( \frac{\Delta}{\Delta + 1} \right) n_i,
\]

\noindent
noting that the vertex $v_i$ is an $S_i$-forcing vertex that forces the vertex $w_i$ to be colored in the first step of the forcing process. After the vertex $v_i$ is played in the forcing process, all vertices in $V(G_i)$ become colored. Thus, $S_i$ is a TF-set in $G_i$ satisfying Inequality~(\ref{Eq2}). This completes our discussion of Case~2.

\medskip
We now return to the proof of Theorem~\ref{t:upperbd}. By the way in which the set $S_i$ is constructed for each $i \in [r]$ (see Case~1 and Case~2 above), either $N[v_i] \subseteq S_i$ or $N[v_i] \setminus \{w_i\} \subseteq S_i$ for some neighbor $w_i$ of $v_i$. We call such a vertex $w_i$ an $S_i$-uncolored neighbor of $v_i$. Further, each such set $S_i$ is a TF-set of $G_i$ satisfying Inequality~(\ref{Eq1}) or Inequality~(\ref{Eq2}). We now let
\[
S' = \bigcup_{i=1}^r S_i.
\]
Further if $r < k$, we let
\[
P' = \bigcup_{i=r+1}^k \{v_i\} \hspace{0.5cm} \mbox{and}  \hspace{0.5cm} A' = \bigcup_{i=r+1}^k \{w_i\},
\]

\noindent
and so $A' = N(P')$ and $|A'| = |P'| = k - r$. Let $A''$ be the set of vertices in $A'$ that are adjacent to no vertex of $A' \cup S'$. We note that if $w \in A''$, then since $G$ is connected, the vertex $w$ has a neighbor $w'$ different from its (unique) neighbor in $P'$. Moreover, such a neighbor $w'$ necessarily belongs to the set $A \setminus (A' \cup S')$. We note further that such a vertex $w$ may have many neighbors $w'$ different from its (unique) neighbor in $P'$. However, each such neighbor $w'$ of $w$ is an  $S_i$-uncolored neighbor of $v_i$ for some $i \in [r]$. Thus,
\[
N(A'') \setminus P'  \subseteq A \setminus (A' \cup S').
\]

\noindent
Let $C$ be a minimum set of vertices in $A \setminus (A' \cup S')$ that dominates the set $A''$. By the minimality of the set $C$, each vertex in $C$ dominates a vertex in $A''$ that is not dominated by the other vertices in $C$. For each vertex $x \in C$, let $x''$ be a vertex in $A''$ that is dominated by $x$ but by no vertex of $C \setminus \{x\}$. Further, let
\[
A''' = \bigcup_{x \in C} \{x''\}.
\]
We note that $A''' \subseteq A''$ and $|C| = |A'''|$. Further, the set $C$ dominates $A''$, and so each vertex of $A'' \setminus A'''$ has a neighbor in $C$. Let
\[
S'' = C \cup (A' \setminus A''')
\]
and consider the set
\[
S = S' \cup S''.
\]

\noindent
We show that the set $S$ is a forcing set in $G$. By construction, the graph $G[S]$ contains no isolated vertex.  As shown earlier, each set $S_i$ is a TF-set of $G_i$ satisfying Inequality~(\ref{Eq1})  or Inequality~(\ref{Eq2}) for all $i \in [r]$. As observed earlier, no vertex of $S_i$ is adjacent to any vertex of $A''$ for all $i \in [r]$. In particular, no vertex of $S_i$ is adjacent to any vertex of $A'''$ for all $i \in [r]$. By the way in which the set $S_i$ is constructed for each $i \in [k]$, either $N[v_i] \subseteq S_i$ or $N[v_i] \cap S_i =  N[v_i] \setminus \{w_i\}$ for some neighbor $w_i$  of $v_i$. In the latter case, we call $w_i$ the $S$-uncolored neighbor of $v_i$. Let $Q = \{v_1,\ldots,v_r\}$, and note that if $r = k$, then $Q = P$, while if $r < k$, then $Q \subset P$.

In the first stage of the forcing process, we color all vertices in $N[Q]$ as follows. If $v_1$ has an uncolored $S$-neighbor (namely, $w_1$), then we play the vertex $v_1$ as our first vertex in the forcing process and consider next the vertex $v_2$. Otherwise, if $N[v_1] \subseteq S_1$, then we immediately consider the vertex $v_2$. If $v_2$ has an uncolored $S$-neighbor (namely, $w_2$), then we  play the vertex $v_2$ in the forcing process and consider next the vertex $v_3$. Otherwise, if $N[v_2] \subseteq S_2$, then we immediately consider the vertex $v_3$.  Proceeding in this way, we first play the vertices in  $Q$ in order whenever such a vertex of $P$ has an $S$-uncolored neighbor. In this way, all vertices in $N[Q]$ are colored.

In the second stage of the forcing process we color all $S$-uncolored vertices of $B_1$. Thereafter, we color all $S$-uncolored vertices of $B_2$, and so on, until we finally color all $S$-uncolored vertices of $B_k$. For this purpose for $i \in [r]$, let $s_i$ be a sequence of played vertices, excluding the vertex $v_i$ (which is already played in the first stage or is never played), in the forcing process in the graph $G_i$ that results in all of $V(G_i)$ being colored, as defined in Case~1.1, Case~1.2 or Case~2. Let $s$ be the sequence consisting of the sequence $s_1$, followed by the sequence $s_2$, followed by the sequence $s_3$, and so on, finishing with the sequence $s_r$.
Recall that after the first stage of the forcing process, all vertices in $N[Q]$ are colored. Playing the vertices of the sequence $s_1$ (in order) in the graph $G$ forces all vertices of $V(G_1)$ to be colored in $G$. Thereafter, playing the vertices in the sequence $s_2$ (in order) forces all vertices of $V(G_2)$ to be colored in $G$. Continuing in this way, once the vertices in the sequence $s_{i-1}$ have been played, we play the vertices in the sequence $s_i$ (in order) thereby forcing all vertices of $V(G_i)$ to be colored in $G$ for all $i \in [r] \setminus \{1\}$. We note that this is always possible since no vertex of $S_i$ is adjacent to any vertex of $A'''$ for all $i \in [r]$. This results in all vertices in $V(G_i)$ colored for all $i \in [r]$.

In the third stage of the forcing process, we play each of the vertices of $C$ in turn, thereby coloring all vertices in $A'''$. Once the sequence of vertices in $C$ are played, we play in the fourth and final stage of the forcing process each vertex in the set $A'$ in turn, thereby coloring its (unique) neighbor in $P'$. In this way, starting with the initial set $S$, there is a sequence of vertices that we can play in the forcing process that results in all vertices of $V(G)$ colored. Thus, $S$ is a TF-set of $G$.

It remains for us to verify that the set $S$ satisfies $|S| \le ( \frac{\Delta}{\Delta + 1} ) n$. Recall that $V(G_i) = \{v_i,w_i\}$ for all $i$ where $r+1 \le i \le k$, and so $n_i = 2$ for such values of $i$. By our earlier observations, $|S_i| \le ( \frac{\Delta}{\Delta + 1} ) n_i$ for $i \in [r]$, implying that
\[
|S'| = \sum_{i=1}^r |S_i|  \le \sum_{i=1}^r \left( \frac{\Delta}{\Delta + 1} \right)  n_i.
\]
Moreover, if $r < k$, then since $|C| = |A'''|$ and $|A'| = k - r$, we note that
\[
\begin{array}{lcl}
|S''| & =&  |C| + |A'| - |A'''| \1 \\
& = & k-r  \\
& = & \displaystyle{ \sum_{i=r+1}^k 1 } \1 \\
& = & \displaystyle{  \sum_{i=r+1}^k \frac{1}{2} n_i } \\
& < & \displaystyle{ \sum_{i=r+1}^k \left( \frac{\Delta}{\Delta + 1}\right) n_i.}
\end{array}
\]

Recall that $n = \displaystyle{ \sum_{i=1}^k n_i}$. We note that either $r = k$, in which case
\[
|S| = |S'| = \sum_{i=1}^k |S_i|  \le \sum_{i=1}^k \left( \frac{\Delta}{\Delta + 1} \right)  n_i = \left( \frac{\Delta}{\Delta + 1} \right) n
\]
or $r < k$, in which case
\[
\begin{array}{lcl}
|S| & = & |S'| + |S''| \1 \\
& < & \displaystyle{  \left( \sum_{i=1}^r \left( \frac{\Delta}{\Delta + 1}\right)  n_i  \right)  + \left( \sum_{i=r+1}^k \left( \frac{\Delta}{\Delta + 1}\right) n_i \right) } \2 \\
& = & \displaystyle{  \sum_{i=1}^k \left( \frac{\Delta}{\Delta + 1}\right)  n_i  } \2 \\
& = & \displaystyle{ \left( \frac{\Delta}{\Delta + 1} \right) n.}
\end{array}
\]

Thus, $F_t(G) \le |S| \le ( \frac{\Delta}{\Delta + 1} ) n$, as desired.
Suppose next that $F_t(G) = ( \frac{\Delta}{\Delta + 1} ) n$. Thus, $S$ is a minimum TF-set in $G$, and $|S| = ( \frac{\Delta}{\Delta + 1} ) n$.
Recall that by our earlier assumptions, $\Delta \ge 3$. %
If $r < k$, then as shown above $|S| < ( \frac{\Delta}{\Delta + 1} ) n$, a contradiction. Hence, $r = k$, implying that $|S_i| = ( \frac{\Delta}{\Delta + 1} ) n_i$ for all $i \in [r]$. This implies that the set $S_i$ must have been constructed in Case~2 for all $i \in [k]$; that is, $B_i = \emptyset$ for all $i \in [k]$. Thus, $B = \cup_{i=1}^k B_i = \emptyset$, implying that the packing $P$ is a perfect packing. By our choice of the packing $P$, this implies that every maximum packing in $G$ is a perfect packing.

Suppose that $v_i$ is a vertex in $P$ of degree at least~$2$ for some $i \in [k]$. If the set $S_i$ is constructed as in Case~1.2, then strict inequality holds in Inequality~(\ref{Eq2}) since in this case $\alpha_i > 1$ implying that $\alpha_i + \ell_i + 1 < \Delta \alpha_i$. Hence, the set $S_i$ must have been constructed as in Case~2. Further  for equality to hold in Inequality~(\ref{Eq2}), the vertex $v_i$ has maximum possible degree, namely~$\Delta$; that is, $d_G(v_i) = \alpha_i = \Delta$. By the way in which $S_i$ is constructed (see Case~2), the set $S_i$ contains all but one vertex from the sets $N[v_i]$, and so $|S_i| = n_i - 1 = \Delta$.

Suppose that $d_G(v_i) = 1$ and $B_i \ne \emptyset$ for some $i \in [k]$. Adopting our earlier notation, this implies that the (unique) neighbor $w_i$ of $v_i$ has degree~$\Delta$ and every neighbor of $w_i$ in $G$ different from $v_i$ belongs to $B_i$. Renaming the vertices in $P$ if necessary, we may choose $i = k$. If $k = 1$, then we obtain a contradiction since the packing consisting of the singelton  vertex $w_1$ of degree~$\Delta$ would contradict our choice of $P$ which currently consists of the singelton  vertex $v_1$ of degree~$1$. Hence, $k \ge 2$.

Since $G$ is connected, we can choose the neighbor $w_k'$ of $w_k$ so that it has at least one neighbor outside $V(G_k)$. Since $w_k' \in B_k$, we note that $w_k'$ is not adjacent with a neighbor of $v_j$ for any $j < k$ for otherwise $w_k' \in B_j$ for some $j \in [k-1]$, a contradiction. However, if remove the vertex $v_k$ from the set $S_k = V(G_k) \setminus \{w_k'\}$ and immediately before playing the vertex $w_k$ in the forcing process we play a neighbor of $w_k'$ outside $V(G_k)$ (which forces $w_k'$ to be colored) and then play the vertex $w_k$ (which forces $v_k$ to be colored), we produce a new TF-set of $G$ of cardinality~$|S| - 1$, a contradiction.

Hence, $d_G(v_i) \ge 2$ for all $i \in [k]$. Suppose that $k \ge 2$. (Recall that $|P| = k$.) We note that $(N[v_1],N[v_2], \ldots, N[v_k])$ is a partition of $V(G)$. As observed earlier, each vertex $v_i$ has~$\Delta$ neighbors for all $i \in [k]$; that is, $|N(v_i)| = \Delta$. Since $G$ is connected, there is a neighbor, $u_1$ say, of $v_1$ that is adjacent to a vertex outside $N[v_1]$. Renaming indices if necessary, we may assume that $u_1$ and $v_2$ have a common neighbor, say $w_2$. Since $u_1$ has degree at most~$\Delta$ and is adjacent to $v_1$, at least one neighbor of $v_i$, say $u_i$, is not adjacent to $u_1$ for all $i \in [k] \setminus \{1\}$.

We now consider the set $T = V(G) \setminus \{u_1,u_2,\ldots,u_k,w_2\}$. Since $\Delta \ge 3$, the graph $G[T]$ is isolate-free noting that $v_j$ and at least one of its neighbors belong to $T$ for all $j \in [k]$. We show that $T$ is a forcing set. The first vertex played in the forcing process is the vertex $v_1$ which colors $u_1$. The second vertex played is the vertex $u_1$ which colors $w_2$. As this stage of the forcing  process all vertices are colored except for the vertices $u_2, \ldots, u_k$. For $i \in [k] \setminus \{1\}$, we play as the $(i+1)$st move in the forcing process the vertex $v_i$ which colors the vertex $u_i$. Once the vertices $v_2, \ldots, v_k$ have been played (after $k+1$ steps in the forcing process)  all vertices of $V(G)$ are  colored. Hence, $T$ is TF-set of $G$. However, $|T| = n - (k+1) = |S| - 1$, and so $|T| < |S| = F_t(G)$,  a contradiction. Therefore, $k = 1$.

Hence, $k = 1$ and the vertex $v_1$ has degree~$\Delta$. If $G \ncong K_{\Delta + 1}$ and $G \ncong K_{1,\Delta}$, then $v_1$ has three neighbors $x_1$, $y_1$ and $z_1$ such that $x_1$ is adjacent to $y_1$ but not to $z_1$. In this case, the set $S^* = V(G) \setminus \{y_1,z_1\}$ is a TF-set of $G$, noting that starting with the set $S^*$ we play as our first vertex in the forcing process the vertex $x_1$ which forces the vertex $y_1$ to be colored, and as our second vertex the vertex $v_1$ which forces the vertex $z_1$ to be colored. Thus, $F_t(G) \le |S^*| < |S| = F_t(G)$, a contradiction. Therefore, $G \cong K_{\Delta + 1}$ or $G \cong K_{1,\Delta}$. This completes the proof of Theorem~\ref{t:upperbd}.~\qed


\medskip
\begin{thebibliography}{99}

\bibitem{AIM-Workshop} AIM Special Work Group, Zero forcing sets and the minimum rank of graphs. \textit{Linear Algebra Appl.} \textbf{428}(7) (2008), 1628--1648.

\bibitem{k-Forcing} D. Amos, Y. Caro, R. Davila, and R. Pepper, Upper bounds on the $k$-forcing number of a graph. \textit{Discrete Applied Math.} \textbf{181} (2015), 1--10.

\bibitem{Barioli13} F. Barioli,  W. Barrett, S. M. Fallat, T. Hall, L. Hogben, B. Shader, P. van den Driessche, and H. van der Holst, Parameters related to tree-width, zero forcing, and maximum nullity of a graph. \textit{J. Graph Theory} \textbf{72}(2) (2013), 146--177.

\bibitem{quantum1} D. Burgarth and V. Giovannetti, Full control by locally induced relaxation. \textit{Physical Review Letters} \textbf{99}(10) (2007), 100501.

\bibitem{logic1} D. Burgarth, V. Giovannetti,  L. Hogben,  S. Severini,  and M. Young. Logic circuits from zero forcing, manuscript. arXiv preprint arXiv:1106.4403, 2011.

\bibitem{CF Complexity} B. Brimkov, Complexity and computation of connected zero forcing, manuscript. arXiv preprint  arXiv:1611.02379, 2016.

\bibitem{Brimkov Davila} B. Brimkov and R. Davila, Characterizations of the connected forcing number of a graph, manuscript. arXiv preprint arXiv:1604.00740, 2016.

\bibitem{Dynamic Forcing} Y. Caro and R. Pepper, Dynamic approach to k-forcing. {\em Theory and Applications of Graphs}, Volume 2: Iss. 2, Article 2, 2015.

\bibitem{ChSu05} S. Chandran and C. Subramanian, Girth and treewidth. \textit{J. Combin. Theory B} \textbf{93} (2005), 23--32.

\bibitem{zf_np} C. Chekuri and N. Korula, A graph reduction step preserving element-connectivity and applications. {\em Automata, Languages and Programming}, 254--265. Springer 2009.


\bibitem{Davila} R. Davila, Bounding the forcing number of a graph. {\em Rice University Masters Thesis}, 2015.

\bibitem{Davila Henning} R. Davila and M. A. Henning, The forcing number of graphs with a given girth, manuscript. arXiv preprint arXiv:1605.02124, 2016.

\bibitem{DaHeMaPe16} R. Davila, M. A. Henning, C. Magnant, and R. Pepper,  Bounds on the connected forcing number of a graph, manuscript. arXiv preprint arXiv:1605.02124, 2016.

\bibitem{Davila Kenter} R. Davila and F. Kenter, Bounds for the zero forcing number of a graph with large girth. {\em Theory and Applications of Graphs}, Volume 2, Issue 2, Article 1, 2015.

\bibitem{girthTheorem} R. Davila, T. Malinowski, and S. Stephen, Proof of a conjecture of Davila and Kenter regarding a lower bound
for the forcing number in terms of girth and minimum degree, manuscript 2016.

\bibitem{GaKwLa93} D. K. Garnick, Y. H. Harris Kwong, and F. Lazebnik, Extremal graphs without three-cycles or four-cycles. \textit{J. Graph Theory} \textbf{17}(5) (1993), 633--645.

\bibitem{Genter1} M. Gentner, L. D. Penso, D. Rautenbach, and U. S. Souzab, Extremal values and bounds for the zero forcing number. \textit{Discrete Applied Math.} \textbf{214} (2016), 196--200.

\bibitem{Genter2} M. Gentner and D. Rautenbach, Some bounds on the zero forcing number of a graph, manuscript. arXiv preprint  arXiv:1608.00747, 2016.

\bibitem{Edholm} C. Edholm, L. Hogben, J. LaGrange, and D. Row, Vertex and edge spread of zero forcing number, maximum nullity, and minimum rank of a graph. \textit{Linear Algebra Appl.} \textbf{436}(12) (2012), 4352--4372.

\bibitem{DomBook1} T. W. Haynes, S. T. Hedetniemi and P. J. Slater, Fundamentals of Domination in Graphs. Marcel Decker, Inc., NY, 1998.

\bibitem{DomBook2} T. W. Haynes, S. T. Hedetniemi and P. J. Slater, Domination in Graphs: Advanced Topics. Marcel Decker, Inc., NY, 1998.

\bibitem{powerdom3} T. W. Haynes,  S. T. Hedetniemi, S. T. Hedetniemi, and M. A. Henning, Domination in graphs applied to electric power networks. \textit{SIAM J. Discrete Math.} \textbf{15}(4) (2002), 519--529.

\bibitem{He98} M. A. Henning, Distance domination in graphs. \emph{Domination in Graphs: Advanced Topics}, T. W. Haynes, S. T. Hedetniemi, and P. J. Slater (eds), Marcel Dekker, Inc. New York, 1998, 335--365.

\bibitem{MHAYbookTD} M. A. Henning and A. Yeo, \emph{Total domination in graphs (Springer Monographs in Mathematics)}.  ISBN-13: 978-1461465249, 2013.

\bibitem{Hogben16} L. Hogbena, M. Huynh, N. Kingsley, S. Meyer, S. Walker, and M. Young, Propagation time for zero forcing on a graph. \textit{Discrete Applied Math.} \textbf{160} (2012), 1994--2005.

\bibitem{LuTang} L. Lu, B. Wu, and Z. Tang, Note: Proof of a conjecture on the zero forcing number of a graph. \textit{Discrete Applied Math.} \textbf{213} (2016), 233--237.



\bibitem{Ore's Theorem} O. Ore, Theory of graphs, Amer. Math. Spc. Colloq. Publ. 3 (1962).

\bibitem{zf_np2} M. Trefois and J. C. Delvenne. Zero forcing sets, constrained matchings and minimum rank. \textit{Linear and Multilinear Algebra}, 2013.

\bibitem{powerdom2} M. Zhao,  L. Kang, and G. Chang, Power domination in graphs. \textit{Discrete Math.} \textbf{306} (2006), 1812--1816.


\end{thebibliography}
\end{document}